\documentclass[12pt,a4paper]{amsart} 

\ifx\pdftexversion\undefined
\usepackage[dvips]{graphicx,color}
\else
\usepackage[pdftex]{graphicx,color}
\fi

\usepackage{amssymb}

\let\ifanglais\iftrue

\def\R{{\mathbb R}}
\def\N{{\mathbb N}}
\def\C{{\mathbb C}}

\def\pP{{\mathbb P}}
\def\vol{{\rm Vol}}

\newcommand{\n}[1]{\Vert #1 \Vert}

\newcommand{\ed}[1]{\textrm{d} #1}
\renewcommand{\leq}{\leqslant} \renewcommand{\geq}{\geqslant}

\newtheoremstyle{mesthm}
  {10pt plus 1pt minus 1pt}
  {9pt minus 6pt}
  {\slshape}
  {0.5cm}
  {\bfseries}
  {.}
  {1ex}
  {}
\newtheoremstyle{mesdefi}
  {6pt plus 1pt minus 1pt}
  {6pt plus 1pt minus 1pt}
  {}
  {0.5cm}
  {\bfseries}
  {.}
  {1ex}
  {}%

\theoremstyle{mesthm}
\newtheorem{lema}{\ifanglais{\large L}emma\else{\large L}emme\fi}
\newtheorem{theo}[lema]{\ifanglais{\large T}heorem\else {\large
    T}h\'eor\`eme\fi}

\newtheorem{prop}[lema]{{\large P}roposition}  
\newtheorem{cor}[lema]{{\large C}orollary}
\newtheorem{ppt}[lema]{\ifanglais{\large P}roperty\else{\large
    P}ropri\'et\'e\fi} 

\newtheorem{conj}{\ifanglais{\large C}onjecture\else{\large C}onjecture\fi}
\newtheorem{rmq}[lema]{\ifanglais{\large R}emark\else{\large
    R}emarque\fi}

\theoremstyle{mesdefi}
\newtheorem{defi}[lema]{\ifanglais{\large D}efinition\else{\large
    D}\'efinition\fi} 
\newtheorem{exs}[lema]{\ifanglais{\large E}xamples\else{\large
    E}xemples\fi} \newtheorem{rem}[lema]{\ifanglais{\large
    R}emark\else{\large  R}emarque\fi}

\title[Amenability in Hilbert geometries]{Spectral radius and Amenability in Hilbert Geometries}
\author{Constantin Vernicos}
\address{%
  Constantin Vernicos \\
Institut de Math\'ematiques\\ 
et de Mod\'elisation de Montpellier\\
UMR 5149 CNRS\\
Universit\'e Montpellier II\\
Case Courrier 051\\
Place Eug\`ene Bataillon\\
F-34095 MONTPELLIER Cedex 5\\
France} 
\email{Constantin.Vernicos@math.univ-montp2.fr}

\keywords{Hilbert Geometry, Finsler, Amenability, Cheeger Constant}
\subjclass[2000]{primary 53C60, secondary 53C24,51F99,53A40} 

\begin{document}

\begin{abstract}
We study the bottom of the spectrum in Hilbert geometries,
we show that it is zero if and only if the geometry is
amenable, in other words if and only if it admits a F\"olner sequence.
We also show that the bottom of the spectrum admits an upper bound, which
depends only on the dimension and
which is the bottom of the spectrum of the Hyperbolic geometry of the same
dimension. Horoballs, from a purely metric point of view, and
their relation with the bottom of the spectrum in Hilbert geometries are briefly
studied.
\end{abstract}

\maketitle
\section*{Introduction and statement of results}

For a Riemanniann manifolds of Ricci curvature bounded from below and
positive injectivity radius it is known thanks to the work of P. Buser \cite{buser},
that the bottom of the spectrum and the Cheeger constant are equivalent and thanks
to M. Kanai \cite{kanai} that the manifolds is quasi-isometric to any of its
discretisations, and that positivity of the Cheeger constant of any discretisation is equivalent
to the positivity of the manifold's Cheeger constant.

The aim of this paper is to prove that such results holds in the setting of
Hilbert geometries. 

Before explaining in more details our results let us recall
what are the objects studied here.

A Hilbert geometry
$(\mathcal{C},d_\mathcal{C})$ is a non empty bounded open convex set $\mathcal{C}$
on $\R^n$ (that we shall call \textit{convex domain}) with
the Hilbert distance 
$d_\mathcal{C}$ defined as follows : for any distinct points $p$ and $q$ in $\mathcal{C}$,
the line passing through $p$ and $q$ meets the boundary $\partial \mathcal{C}$ of $\mathcal{C}$
at two points $a$ and $b$, such that one walking on the line goes consecutively by $a$, $p$, $q$
$b$. Then we define
$$
d_{\mathcal C}(p,q) = \frac{1}{2} \ln [a,p,q,b],
$$
where $[a,p,q,b]$ is the cross ratio of $(a,p,q,b)$, i.e., 
$$
[a,p,q,b] = \frac{\| q-a \|}{\| p-a \|} \times \frac{\| p-b \|}{\| q-b\|} > 1,
$$
with $\| \cdot \|$ the canonical euclidean norm in
$\mathbb R^n$.

Note that the invariance of the cross ratio by a projective map implies the invariance 
of $d_{\mathcal C}$ by such a map.

These geometries are naturally endowed with
a  $C^0$ Finsler metric $F_\mathcal{C}$ as follows: 
if $p \in \mathcal C$ and $v \in T_{p}\mathcal C =\R^n$
with $v \neq 0$, the straight line passing by $p$ and directed by 
$v$ meets $\partial \mathcal C$ at two points $p_{\mathcal C}^{+}$ and
$p_{\mathcal C}^{-}$~; we then define
$$
F_{\mathcal C}(p,v) = \frac{1}{2} \| v \| \biggl(\frac{1}{\| p -
  p_{\mathcal C}^{-} \|} + \frac{1}{\| p - p_{\mathcal C}^{+}
  \|}\biggr) \quad \textrm{and} \quad F_{\mathcal C}(p , 0) = 0.
$$ 


The Hilbert distance $d_\mathcal{C}$ is the length distance associated to 
$F_{\mathcal C}$.

Thanks to that Finsler metric, we can built a Borel measure
$\mu_{\mathcal C}$ on $\mathcal C$ (which is actually
the Hausdorff measure of the metric space $(\mathcal C,
d_{\mathcal C})$, see \cite{bbi}, exemple~5.5.13 ) as follows.

To any $p \in \mathcal C$, let $B_{\mathcal C}(p) = \{v \in \R^n ~|~ F_{\mathcal{C}}(p,v) < 1 \}$
be the open unit ball in
$T_{p}\mathcal{C} = \R^n$ of the norm $F_{\mathcal{C}}(p,\cdot)$ and 
$\omega_{n}$ the euclidean volume of the open unit ball of the standard euclidean space
$\R^n$.
Consider the (density) function $h_{\mathcal C}\colon  \mathcal C \longrightarrow \R$ given by $h_{\mathcal C}(p)
= \omega_{n}/\vol\bigl(B_{\mathcal C}(p)\bigr),$ where $\vol$ is the canonical Lebesgue measure
of $\R^n$. We define $\mu_{\mathcal C}$, which we shall call
the \textit{Hilbert Measure} on $\mathcal{C}$,
by
$$
\mu_{\mathcal C}(A) = \int_{A} h_{\mathcal C}(p) \ed{\vol(p)}
$$
for any Borel set $A$ of $\mathcal C$.

The bottom of the spectrum of $\mathcal{C}$, denoted by $\lambda_1(\mathcal C)$,
is defined as in a Riemannian manifold of infinite volume, thanks to
the Raleigh quotients as follows
\begin{equation}
\label{eqdef1}
\lambda_1(\mathcal{C})=
\inf \frac{\displaystyle{\int_\mathcal{C}{\n{df_p}_\mathcal{C}^*}^2~d\mu_\mathcal{C}(p)}}
{\displaystyle{\int_\mathcal{C} f^2(p)d\mu_\mathcal{C}(p)}}\text{,}
\end{equation}
where the infimum is taken over all non zero lipschitz functions with
compact support in $\mathcal{C}$

Finally the Cheeger constant of $\mathcal C$ is defined by
 \begin{equation}
\label{Cheeger}
I_{\infty}(\mathcal{C})=
\inf_U \frac{\nu_{\mathcal C}(\partial U)}{\mu_{\mathcal C}(U)}\text{,}
\end{equation}
where $U$ is an open set in $\mathcal{C}$ whose closure is compact and
whose boundary is a $n-1$ dimensional submanifold, and $\nu_{\mathcal C}$ is the
Hausdorff measure associated to the restriction of the finsler norm $F_\mathcal{C}$
to hypersurfaces.

When the convex set is a euclidean ball one gets the Klein or projective model
of the Hyrerbolic geometry. Hence 
one of the objects of research in Hilbert geometries is to understand how close
they can be to the Hyperbolic geometries. A property which has been studied
a lot recently is gromov hyperbolicity (for related papers on hyperbolicity 
in Hilbert geometries see \cite{benoist}, \cite{benoist1},\cite{kn} and \cite{cvv}).

When we began the study of the spectrum in Hilbert Geometry, we started by looking
at plane Hilbert Geometries in \cite{cvc}. There we found out that the positivity of the bottom
of the spectrum was equivalent to the hyperbolicity in the sens of Gromov. Two main
ingredients were involved. The first one is that in the two dimensional case, if the boundary
of the Hilbert geometry is not strictly convex, then the bottom of the spectrum is zero.
The seconde one consisted in proving the equivalence for the Cheeger constant, and then
thanks to a Cheeger type inequality deduce it for the bottom of the spectrum.

In Higher dimension we finally found out in \cite{cvc2} 
that a Hilbert geometry did not need to be strictly
convex to have a positive bottom of the spectrum. However by showing that
the Hilbert geometries had bounded local geometry and using a paper of J. Cao \cite{cao}\footnote{Actually, J. Cao uses a theorem of M. Kanai to conclude that the positivity
of the Cheeger constant of his space is equivalent the positivity of the cheeger constant of
some discrete metric space to which his space is quasi-isometric. However Kanai's theorem does not apply in J. Cao setting. Hence one should be careful while using his theorem, or one might make a mistake. In the present paper we prove that M. Kanai results holds in the setting of Hilbert Geometries, which fully justifies our result in \cite{cvc2}.},
we were able to prove that if the geometry was hyperbolic in the sens of Gromov, once
again the Cheeger constant had to be positive and by our Cheeger type inequality deduce
the same for the bottom of the spectrum.

There however was a missing link to clarify what makes the bottom of the spectrum
zero. Then one thinks of two types of results. The first one, mentioned at the
beggining, is due to P. Buser \cite{buser}
who shows that in Riemanniann geometry, under the right assumptions on the curvature and injectivity radius, 
the positivity of bottom of the spectrum and that of
the Cheeger constant is in fact equivalent. The second one
is due to the late R. Brooks \cite{brooks} who shows that the bottom of the spectrum of
the covering of a compact Riemannian manifold is positive if and only if its fundamental group
is not amenable.

If amenability makes sense for a divisible Hilbert geometry (which admits
a group of isometry which acts cocompactly on it), in the general case there is no group
big enough to do anything \cite{so}. However for a discrete metric spaces, one may require
the pseudo group of bounded perturbations of the identity to be amenable \cite{cgh} 
(see also section \ref{amenability} in the present paper).
For such metric spaces, similar results combining the equivalence of R. Brooks
and P. Buser exist under suitable conditions \cite{cgh} 

Hence we are naturally led to say that a Hilbert Geometry is amenable if and only
if it is quasi-isometric to a discrete metric space which is amenable. Taking that path
and in the light of
M. Kanai paper \cite{kanai}, we are bound to study discretisations of the Hilbert
Geometry themselves. This led us to our first result

\begin{theo}
  Let $(\mathcal{C},d_\mathcal{C})$ be a Hilbert geometry, then it is quasi-isometric
to any of its discretisation, and thus any two of its discretisations are
quasi-isometric.
\end{theo}

Thanks to this first result we see that focusing on a discretisation is a good idea, for
amenabilty is invariant by quasi-isometry.
Furthermore these discretisations are also of bounded geometry and thus
the classical results linking amenability, spectral radius of a simple random walk
and the cheeger constant apply to them. However we still have to climb back to
the Hilbert geometry. This is possible thanks the local boundedness of the geometry
proved in \cite{cvc2}, and we finally obtain 

\begin{theo}[Main Theorem]
  Let $(\mathcal{C},d_\mathcal{C})$ be a Hilbert geometry then the following are equivalent
  \begin{enumerate}
      \item The bottom of the spectrum of $\mathcal{C}$ si positive;
      \item The spectral radius of any discretisation is less than $1$;
      \item The Cheeger constant of $\mathcal{C}$ is positive;
      \item The Cheeger constant of any discretisation is positive;
      \item $\mathcal{C}$ is not amenable.
  \end{enumerate}
\end{theo}

The strategy consists in showing the equivalence between (1) and (2) (which is the
content of section \ref{lambdadiscret}),
(3) and (4) (we do it in section \ref{isodiscret}) 
and showing that a discretisation has the good property for
(2), (4) and (5) to be equivalent (see  \cite{cgh} and section \ref{isobot} of this paper). For the convenience
of the reader, these equivalences are proved in full details in the setting
of Hilbert Geometries. However one will find out while reading our proof that in fact
the real important property is the fact that Hilbert Geometries are of local bounded geometry.
 
In section \ref{polygones} I also introduce a familly of convex sets whose Hilbert geometry is amenable which I call $G_n$-polygons (where $G_n$ stands for $PGL_n(\R)$). I believe that they are the only ones to have a Hilbert geometry
which is amenable.

After focusing on the lower bound, it's logic to focus on the upper bound. In this
paper we answer to the first part of a question of B.~Colbois in the following way
\begin{theo}[Upper bound of the spectrum]
  Let $(\mathcal{C},d_\mathcal{C})$ be a Hilbert Geometry with $\mathcal{C}$ a bounded open
convex set in $\R^n$. Then 
$$
\lambda_1(\mathcal{C}) \leq \frac{(n-1)^2}{4}\text{.}
$$
\end{theo}

Hence now the second part of the question makes sense: Is there a rigidity involved in that
equality, i.e., is the equality only achieved by the Hyperbolic geometries ?

In the first section of this paper, we also study Horoballs, in a purely metric
point of view (no dynamics, sorry !) and their links with the bottom of the spectrum
(there one can replace it with the Cheeger constant and obtain the same kind of results).

\noindent\textbf{Acknowledgment:} I wish to thank A. Valette for his never-ending patience
in answering my questions on amenability and pointing out to me the paper \cite{cgh}.
\section{Horoballs of Hilbert Geometries}

\begin{quote}
In this section we give a definition of Horoballs, some examples, 
and study their relationship with
the bottom of the spectrum.
\end{quote}

\begin{defi}
  Let $\mathcal{C}$ be a properly open convex set in $\pP^n$. We will
call $\mathcal{H}$ an horoball of $\overline{\mathcal{C}}$ if there exists
a point $x_0$ and $p\in \partial\mathcal{C}$ such that the familly of balls
$B_{\mathcal{C}}\bigl(x,d_\mathcal{C}(x,x_0)\bigr)$ where $x \in (x_0,p)$
converges to $\mathcal{H}$ as $x\to p$ for the hausdorff topology of $\pP^n$. 
We say that $\partial\mathcal{H}$ is the horosphere 
based at $p$ passing through $x_0$. We may some time denote this by $\mathcal{H}_{p,x_0}$.
\end{defi}

\begin{ppt}
  For any points $(x_0,p)\in \mathcal{C}\times \partial \mathcal{C}$, there is a Horosphere
based at $p$ passing by $x_0$.
\end{ppt}
\begin{proof}
  Let $x,x'$ in $(x_0,p)$ such that 
$$
d_{\mathcal{C}}(x,x_0)<d_{\mathcal{C}}(x',x_0)
$$
and let $y\in B_\mathcal{C}\bigl(x,d_\mathcal{C}(x,x_0)\bigr))$,
which means that $d(y,x)\leq d_{\mathcal{C}}(x,x_0)$.
Then
$$
d(y,x')\leq d(y,x) +d(x,x') \leq  d_{\mathcal{C}}(x,x_0)+d_\mathcal{C}(x,x')
$$
as $x$, $x'$ and $x_0$ are on the same line we obtain
$$
d_{\mathcal{C}}(x,x_0)+d_\mathcal{C}(x,x')=d_\mathcal{C}(x',x_0)
$$
thus $y\in B_\mathcal{C}\bigl(x',d_\mathcal{C}(x',x_0)\bigr))$.
Hence the familly of ball $B_{\mathcal{C}}\bigl(x,d_\mathcal{C}(x,x_0)\bigr)$ is
increasing, and bounded, thus it converges to some subset of $\overline{\mathcal{C}}$.
\end{proof}

\begin{exs}
The following figure illustrate the previous proof in a triangle
$$
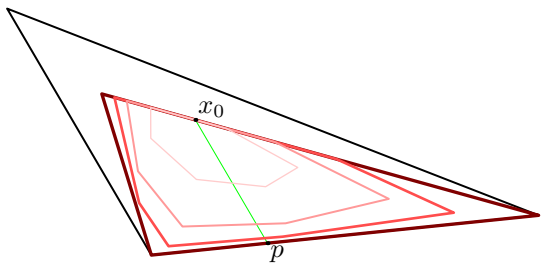
\includegraphics[scale=.5]{horoball.2}
$$
In an hexagone, a Horosphere looks like that:
$$
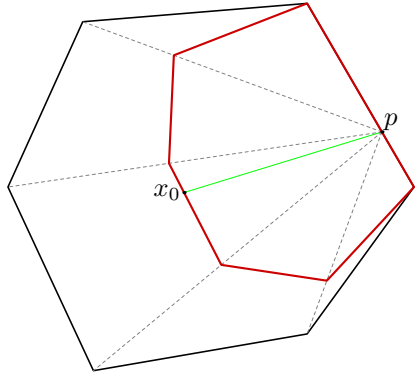
\includegraphics[scale=.4]{horoball.3}
$$
The following gives exemples of Horoballs in $\mathcal{C}=\{x^{4}+y^{1.1}\leq 1 \}$
centered at the same point $p$
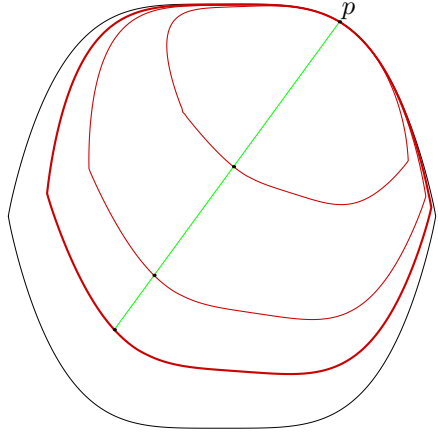
$$\includegraphics[scale=.4]{horoball.4}
$$
 
\end{exs}

\begin{prop}
  For a properly open convex, 
$$
\lambda_1(\mathcal{C})=\inf_{\mathcal{H}}\lambda_1(\mathcal{H})
$$
where the infimum is taken over all horoballs $\mathcal{H}$.
\end{prop}
\begin{proof}
  We just need to prove that
  \begin{equation}
    \label{eqlautre}
    \lambda_1(\mathcal{C})\geq\inf_{\mathcal{H}}\lambda_1(\mathcal{H})\text{.}
  \end{equation}
Let us fix a point $x_0$. Then for any $R$, taking a line passing by $x_0$ it crosses
the boundary of $\partial\mathcal{C}$ at two points $p$ and $q$ and the ball $B_\mathcal{C}(x_0,R)$
at $y$ and $x$. Let us suppose that the point on the line are consecutively $p$,$x$,$x_0$,$y$ and $q$.
Then the ball $B_\mathcal{C}(x_0,R)$ is inside the Horoball based at $p$ passing by $y$.
Hence 
$$
\lambda_1(B_\mathcal{C}(x_0,R))\geq\inf_{\mathcal{H}}\lambda_1(\mathcal{H})\text{.}
$$
Passing to the limit in $R$ we get our result.
\end{proof}

Recall that a convex domain is said to be \textsl{divisible} if their
is a subgroup of isometries acting co-compactly on it.
Let us make the last statement, in the divisble case, more precise. 

\begin{prop}
  For a divisible convex set $\mathcal{C}$, 
let $\mathcal{H}$ be a Horoball then,
$$
\lambda_1(\mathcal{C})=\lambda_1(\mathcal{H})
$$
\end{prop}
\begin{proof}
  What is clear is that 
  \begin{equation}
    \label{eqpremier}
\lambda_1(\mathcal{C})\leq \lambda_1(\mathcal{H})\text{.}    
  \end{equation}

Consider a point $x_0$ in $\mathcal{C}$, and a point $p$ on the boundary.
Let $x$ be a point on the segment $(x_0,p)$ and consider the $d_\mathcal{C}$ ball
centered at $x$ and passing by $x$. As $x\to p$, this balls converges to the
horoball passing by $x_0$ based at $p$.

Furthermore, by cocompactness, if $\Gamma$ is a group which divides $\mathcal{C}$,
then there exists some constant $C$, such that for every $x \in (x_0,p)$,
there exist $\gamma\in\Gamma$ such that
$$
d_\mathcal{C}(x,\gamma x_0)\leq C
$$
hence for any $x$ such that $d_\mathcal{C}(x,x_0)>C$ we have,
$$
B_\mathcal{C}\bigl(x,d_{\mathcal{C}}(x,x_0)-C\bigr)\subset B_\mathcal{C}\bigl(\gamma x_0,d_{\mathcal{C}}(x,x_0)\bigr)\subset B_\mathcal{C}\bigl(x,d_{\mathcal{C}}(x,x_0)+C\bigr)
$$
by which we deduce that
\begin{multline}
  \lambda_1\Bigl(B_{\mathcal{C}}\bigl(x,d_{\mathcal{C}}(x,x_0)+C\bigr)\Bigr)\leq\\
\lambda_1\Bigl(B_{\mathcal{C}}\bigl(x_0,d_{\mathcal{C}}(x,x_0)\bigr)\Bigr)\leq\\ 
\lambda_1\Bigl(B_{\mathcal{C}}\bigl(x,d_{\mathcal{C}}(x,x_0)-C\bigr)\Bigr)
\end{multline}

this implies that
$$
\lim_{x\to p} \lambda_1\Bigl(B_{\mathcal{C}}\bigl(x_0,d_{\mathcal{C}}(x,x_0)\bigr)\Bigr)=\lim_{x\to p} \lambda_1\Bigl(B_{\mathcal{C}}\bigl(x,d_{\mathcal{C}}(x,x_0)\bigr)\Bigr) \text{.}
$$

Now notice that
$$
\lim_{x\to p} \lambda_1\Bigl(B_\mathcal{C}\bigl(x_0,d_{\mathcal{C}}(x,x_0)\bigr)\Bigr) 
= \lim_{R\to\infty }\lambda_1\bigl(B_\mathcal{C}(x_0,R)\bigr) = \lambda_1(\mathcal{C})\text{.}
$$

(This comes the fact that if $f_k$ is a sequence of functions with compact support such
that their Rayleigh quotient converges to the bottom of the spectrum,
then we can find a sequence of balls with increasing radius on which they are defined)

Now let us finally notice that
$$
\lim_{x\to p} \lambda_1\Bigl(B_\mathcal{C}\bigl(x,d_{\mathcal{C}}(x,x_0)\bigr)\Bigr) \geq \lambda_1\bigl(H(x_0,p)\bigr)
$$
From which we deduce, thanks to (\ref{eqpremier}) that
$$
\lambda_1(\mathcal{C})= \lambda_1\bigl(H(x_0,p)\bigr)
$$
which also implies that the right part of this equality neither depends on $p$ nor on $x_0$.
\end{proof}

\section{Upper bound}
\begin{quote}
In this section we give an optimal upper bound on the bottom of the spectrum
of Hilbert geometries.
\end{quote}

Let us recall that by adding a projective hypreplane to $\R^n$, we can see
it as subset of the projective space $\pP^n$. Then, convex domains of $\R^n$ coincide with
open and bounded convex subspaces of $\pP^n$ which doesn't meet at least
one projective hyperplane. We call such subspace of $\pP^n$ \textsl{properly open
convex sets}.

Let us first begin with an easy case :
\begin{lema}
  Let $\mathcal{C}$ be a properly open convex set in $\pP^n$ which 
admits an osculating ellipsoid then 
$$ \lambda_1(\mathcal{C})\leq\dfrac{(n-1)^2}{4}\text{.}$$
\end{lema}

\begin{proof}
  By a result due to Benzecri (see \cite{benzecri} page 325, proposition 10),
if $\mathcal{C}$ admits an osculating ellipsoid $\mathcal{E}$, then there exists
a sequence of projective transformations $g_n\in G_n$ such that $g_n\mathcal{C}$ tends
to $\mathcal{E}$ as $n$ goes to $\infty$. Now a result of Colbois-Vernicos \cite{cvc}, 
implies that the $\lambda_1$ is upper semi-continuous with respect to the haussdorf topology
on properly open convex sets,
hence 
$$
\limsup_{n\to \infty}\lambda_1(g_n\mathcal{C})\leq \lambda_1(\mathcal{E})
$$
however $\lambda_1(\mathcal{E})=(n-1)^2/4$ and for any $n$, $\lambda_1(g_n\mathcal{C})=\lambda_1(\mathcal{C})$,
thus our lemma follows.
\end{proof}

Remark also that the upper semi-continuity implies that the family of convex sets
such that $\lambda_1=(n-1)^2/4$ is not dense in the familly of properly open convexe. 
More precisely, the only familly which is dense, is the familly with
zero $\lambda_1$.

There remains the general case, for this one we will use Alexandroff's theorem
which states that any convex set is almost everywhere two times differentiable.
This implies that at almost every point of the boundary there is a an ellipsoid tangent
and localy inside the convexe (see also Berck-Vernicos \cite{berckver}).

If one considers a point $x_0$ inside the convex and the asymptotic balls, these are the
images of the boundary under the dilation centered at $x_0$ of ratio $\tanh(R)$, and pull
back the finsler area of this asymptotic balls divided by $e^{(n-1)R}$ as $R\to +\infty$ on the boundary,
one gets a measure on the boundary which is in $L^1$ (see  Berck-Vernicos \cite{berckver}).
From Egoroff theorem, this implies, that on the boundary, for any $\eta>0$ there is a set of measure
$\eta$ on the complement of which there is uniforme convergence of these measures to the limit
measure.

This gives the following

\begin{prop}
  Let $\mathcal{C}$ be a properly open convex set in $\pP^n$ then 
$$ \lambda_1(\mathcal{C})\leq\dfrac{(n-1)^2}{4}\text{.}$$
\end{prop}

\begin{proof}
Let $0$ be a fixed point of our convex, and let us denote for any $\varepsilon>0$
$$
h_\varepsilon = (n-1)+\varepsilon \text{.}
$$

Let $\eta>0$ and consider the subset of the boundary $S_\eta$ on the complement 
of which there is uniform convergence
of the density of area of asymptotic spheres.
Let us denote by $B_\eta$ the set of lines from $0$ to the complement of $S_\eta$
on the boundary.

Thanks to the coarea formula and Egoroff's theorem, 
$\exp(-(n-1)\cdot R) {\rm Vol}(B_H(R)\cap B_\eta)$ converges to some number (see also
Berck-Vernicos \cite{berckver}). Hence
$$
\exp(-h_\varepsilon\cdot R) {\rm Vol}(B_H(R)\cap B_\eta)\to 0\text{.}
$$

The idea is a classical one and consists in showing that for $\varepsilon>0$ we have
$$
\lambda_1(\mathcal{C})\leq \frac{h_\varepsilon^2}{4}
$$
passing to the limit our results will then follow.

We will consider the familly of funtions $(F_R)$ defined as follows:
$$
x\mapsto \exp\bigl(- h_\varepsilon/2 d_H(0,x)\bigr) - \exp \bigl(-h_\varepsilon/2 R\bigr) = \varphi_{h_\varepsilon}(x)- \exp (-h_\varepsilon/2 R) 
$$   
on the ball of radius $R$ centered at $0$ and equal to $0$ outside this ball.

We will compute the rayleigh quotient on $B_\eta$,
because we have $\lambda_1(\mathcal{C})\leq \lambda_1^{\mathcal{C}}(B_\eta)$.

Let us compute the differential of the function $F_R$ (we write $\psi(x)=d_H(0,x)$)
where it is not $0$
$$
d_{\mid x}F_R \cdot v = -\frac{h_\varepsilon \varphi_{h_\varepsilon}(x)}{2}d_{\mid x }\psi\cdot v
$$ 

From this we deduce the following expression of the Rayleigh quotient of $F_R$
$$
\mathcal{R}(F_R) =  \frac{h_\varepsilon^2}{4}\dfrac{\int_{B_\mathcal{C}(R)\cap B_\eta} \varphi_{h_\varepsilon}^2(x) d\mu(x)}%
{\int_{B_\mathcal{C}(R)\cap B_\eta} \bigl(\varphi_{h_\varepsilon}(x) -\exp(-h_\varepsilon/2 R)\bigr)^2 d\mu(x)}
$$
Thus to obtain our results it suffices to show that this quotient tends to $h_\varepsilon^2/4$
for a suitable subfamilly of real numbers.

In other words we must show
that as $R\to \infty$ the following quotient tends to something smaller than $1$: 
\begin{multline}
\dfrac{\int_{B_\mathcal{C}(R)\cap B_\eta} \varphi_{h_\varepsilon}^2(x) d\mu(x)}%
{\int_{B_\mathcal{C}(R)\cap B_\eta} \varphi_{h_\varepsilon}^2(x) -2\varphi_{h_\varepsilon}(x)\exp(-h_\varepsilon/2 R) + \exp(-h_\varepsilon\cdot R) d\mu(x)}\\
\leq \dfrac{K(R)}%
{K(R)+P(R)- P(R)^{1/2}K(R)^{1/2}}  
\end{multline}

where 
\begin{eqnarray}
  \label{eqketp}
  K(R) &=& \int_{B_\mathcal{C}(R)\cap B_\eta} \varphi_{h_\varepsilon}^2(x)\,d\mu(x) \text{, and}\\
  P(R) &=&\exp(-h_\varepsilon\cdot R) {\rm Vol}(B_\mathcal{C}(R)\cap B_\eta)
\end{eqnarray}

Hence it suffices to show that $K(R)$ does not go to $0$ as $R\to \infty$, while $P(R)$ does.

By definition of the $B_\eta$ 
we have $\lim_{R \to \infty}\ln(P(R))/R \leq -\varepsilon $. 
This means that for any sequence of real number $(R_n)$ which goes to infinity,
and for some $N\in \N$, then for any $n>N$  we
have $$\ln(P(R_n))/R_n\leq -\varepsilon/2 \text{,}$$ hence $P(R_n)\leq \exp(-\varepsilon R_n/2))$.
As for $K(R_n)$ we have $K(R_n)>K(1)$.

Thus as $n \to \infty$ we deduce that
$$
\lambda_1(\mathcal{C})\leq \lambda_1^{\mathcal{C}}(B_\eta)\leq\frac{h_\varepsilon^2}{4}\text{.}
$$

This being true for any $\varepsilon$ we finally get
$$
\lambda_1(\mathcal{C})\leq\frac{(n-1)^2}{4}\text{.}
$$
\end{proof}

\section{Positivity of the bottom of the spectrum}\label{polygones}

\subsection{$G_n$-Polygons}

\begin{defi}
  Let $\mathcal{C}$ be a properly open convex set in $\pP^n$,  we will say that
 $\mathcal{C}$ is $G_n$-polygonal, if there is a polygone in $\overline{G_n\mathcal{C}}$
\end{defi}

\begin{rem}
  In the two dimensional case, one can replace "polygon" by "triangle", and then by
Y.~Benoist's result in \cite{benoist}
a plane convex set is not $G_n$-polygonal if and only if it is $\delta$-hyperbolic.
\end{rem}

\begin{prop}
  Let $(\mathcal{C},d_\mathcal{C})$ be a $G_n$-polygonal properly open convex set in
$\pP^n$, then the bottom of its spectrum is zero.
\end{prop}

In fact this proposition follows from the following property and the semicontinuity
of the bottom of the spectrum with respect to the Hausdorff topology on compact sets of $\pP^n$.

\begin{ppt}
  The bottom of the spectrum of a polygon is zero.
\end{ppt}
\begin{proof}
  To do this we show that a polygon has polynomial volume growth.
We do this by induction. 
\textbf{Claim:} This is true for $2$-dimensional 
polygones. 

\medskip

Suppose that a all $n$-dimensional polygons have polynomial volume growth and
consider $P_{n+1}$ a $n+1$ dimensional polygons. Choose a point $x_0$ inside $P_{n+1}$.
Now a non trivial argument used and proved in \cite{berckver} (see also \cite{cvp}
and lemma \ref{crofton} in the present paper), says
that the spheres of radius $R$ centered at $x_0$ and the asymptotics spheres obtained
by a dilation  of ratio $\tanh R$ centered at $x_0$ of $P_{n+1}$ have the
same asymptotic behaviour in terms of $n$-volume. However the asymptotic volume
of a face of the asymptotic sphere of ratio $\tanh R$ is of order $R^n$.
This implies the existence of two  constants $C(P_{n+1})_i$, $i=1$ and $2$, such that the $n$-volume
of the sphere of large radii is between $C(P_{n+1})_1\cdot R^n$ and $C(P_{n+1})_2\cdot R^n)$. 
Now using
the co-area inequality showd in Berck-Vernicos \cite{berckver} (see also \cite{cvc}),
one gets that the asympotic volume of the balls of radius $R$  is polynomial of
order $n+1$, i.e., there exists two constants $A$ and $B$ such that
$$
A\cdot R^{n+1}\leq \vol_{\mathcal{C}}\bigl(B(x_0,R)\bigr)\leq B\cdot R^{n+1}\text{.}
$$

Now let us show the claim. This is done by showing that taking a point $x_0$ in $P_2$,
and the asymptotic sphere of radius $R$, then its edges have length
asympotically equal to $2\cdot R$ (easy computation left to the reader).
Hence the asympotic length of a ball is of
order $2$ times the number of sided of $P_2$ times $R$.
Again the co-areas inequality implies that the asymptotic volume of
$P_2$ is of order $R^2$.

Now taking adapted test functions on the balls one easily shows that 
the $\lambda_1$ of our polygons is zero.
\end{proof}

Another consequence is the following.

\begin{prop}
  Let $F_\lambda=\{ \mathcal{C}\in \R^n \mid \lambda_1(\mathcal{C})\geq\lambda\}$, then $F_\lambda$ is a closed, 
$G_n$-invariant such that none of its elements are $G_n$-polygonal.
\end{prop}

\begin{proof}
The only real difficulty lies in the closenesness.
Indeed, let $\mathcal{C}$ be in $F_\lambda$, then it can't be a polygon, 
by the previous property.
Moreover the upper semi-continuity of the $\lambda_1$ (See Colbois-Vernicos \cite{cvc})
implies that for any sequence $\mathcal{C}_n$ in $F_\lambda$
converging to some convex set $\mathcal{C}$ one has
$$
\lambda\leq\limsup_{n\to \infty}\lambda_1(\mathcal{C}_n)\leq\lambda_1(\mathcal{C})
$$
thus $\mathcal{C}$ is in $F_\lambda$, hence $F_\lambda$ is closed.
\end{proof}


\subsection{Amenability}\label{amenability}

We recall some definitions from \cite{cgh}.

\begin{defi}[Pseudogroup of transformation]
  A pseudo group $\mathcal{G}$ of transformations of a set $X$, also
denoted by $(\mathcal{G},X)$ is a set
of bijections $\gamma\colon S\to T$ between subsets $S$, $T$ of $X$ which satisfies
the following conditions
\begin{enumerate}
\item The identity $X\to X$ is in $\mathcal{G}$;
\item if $\gamma\colon S\to T$ is in $\mathcal{G}$, so is the inverse $\gamma^{-1}\colon T\to S$;
\item if $\gamma\colon S\to T$ and $\delta\colon T\to U$ are in $\mathcal{G}$ so is $\delta\circ\gamma\colon S\to U$;
\item if $\gamma\colon S\to T$ is in $\mathcal{G}$ and if $S'\subset S$, the restriction $\gamma_{S'}\colon S'\to \gamma(S')$
is in $\mathcal{G}$;
\item if $\gamma\colon S\to T$ is a bijection between two subsets $S$,$T$ of $X$ and
if there is a \textit{finite} partition of $S=\sqcup_{1\leq j\leq n}S_j$ ($\sqcup$ stands for disjoint union) 
with $\gamma_{S_j}$ in $\mathcal{G}$ for $j\in \{1,\ldots,n\}$ then $\gamma$ is in $\mathcal{G}$. 
\end{enumerate}
For $\gamma\colon S\to T$ in $\mathcal{G}$ we write also $\alpha(\gamma)$ for the domain $S$ of $\gamma$
and $\omega(\gamma)$ for its range $T$.
\end{defi}

In the following definition we denote by $\mathcal{P}(X)$ the set of all subsets of $X$.

\begin{defi}[$\mathcal{G}$-invariant mean]
  A $\mathcal{G}$-invariant mean on $X$ is a mapping $\mu\colon \mathcal{P}(X)\to [0,1]$ which
is
\begin{enumerate}
\item Finitely additive: $\mu(S_1\cup S_2)=\mu(S_1)+\mu(S_2)$ for $S_1$, $S_2\in \mathcal{P}(X)$ with
$S_1\cap S_2=\emptyset$;
\item Invariant: $\mu\bigl(\omega(\gamma)\bigr)=\mu\bigl(\alpha(\gamma)\bigr)$;
\item normalised: $\mu(X)=1$
\end{enumerate}
\end{defi}

Hence we say that the pseudogroup $\mathcal{G}$ is amenable if there
exists a $\mathcal{G}$-invariant mean on $X$.

We are going to focus on a specific pseudogroup associated to metric spaces.

\begin{defi}[The bounded perturbations]
For a metric spave $(X,d)$, the pseudogroup  
$\mathcal{W}(X)$ of \textit{bounded perturbation of the identity} consists
of bijection $\gamma\colon S\to T$, where $S$ and $T$ are subsets of $X$ and
$$
\sup_{x\in S }d\bigl(\gamma(x),x\bigr)< \infty\text{.}
$$  
\end{defi}

Let us recall that a subset $X\in (\mathcal{C},d_\mathcal{C})$ is a separated net
if there exists a constant $r>0$ for which the two following properties hold:
\begin{enumerate}
\item $d_\mathcal{C}(x,y)\geq r$ for all $x,y\in X$, $x\neq y$;
\item $X$ is a maximal subset of $\mathcal{C}$ for this property;
\end{enumerate}

Thus we are now able to define our notion of amenability related
to the Hilbert geometries

\begin{defi}[Amenability]
  Let $(\mathcal{C},d_\mathcal{C})$ be Hilbert geometry, we will say that $\mathcal{C}$
is amenable if and only if, for some separated net $X$ of $\mathcal{C}$, the pseudo group of bounded
perturbation of the identity, $\mathcal{W}(X)$ is
amenable.
\end{defi}

Let us now state our main theorem, which will be proved in the next section

\begin{theo}\label{maintheo}
 Let  $(\mathcal{C},d_\mathcal{C})$  be a Hilbert geometry. Then following are equivalent
 \begin{enumerate}
     \item $(\mathcal{C},d_\mathcal{C})$ is amenable;
     \item $\lambda_1(\mathcal{C})=0$;
     \item $I_\infty(\mathcal{C})=0$.
 \end{enumerate}
\end{theo}


This gives a clearer point of view on the nullity of the bottom of the spectrum.The following results is a consequence of
\cite{cvc} and \cite{cvc2} 

\begin{cor}
  Let $(\mathcal{C},d_\mathcal{C})$ be a Hilbert geometry,
if $(\mathcal{C},d_\mathcal{C})$ is $\delta$-hyperbolic, then it is not amenable.
\end{cor}
\begin{proof}
  If $(\mathcal{C},d_\mathcal{C})$ is $\delta$-hyperbolic, then the  bottom
  of its spectrum is positive.
\end{proof}

\begin{rem}
Notice that if there is a set $\Omega\in \overline{G_n\mathcal{C}}$ which is quasi-isometric 
to an amenable group, then $\mathcal{C}$ is amenable.
\end{rem}

\begin{prop}
  A $G_n$-polygon is amenable.
\end{prop}
\begin{proof}
  Indeed we know that a $G_n$-polygon has a it bottom of the spectrum which is null, hence 
thanks to theorem \ref{maintheo} we know that it is amenable.
\end{proof}

\begin{prop}
  Let $F_\lambda=\{ \mathcal{C}\in \R^n \mid \lambda_1(\mathcal{C})>\lambda\}$, then $F_\lambda$ is a closed, 
$G_n$-invariant such that none of its element is amenable. 
\end{prop}

\begin{proof}
Follows from the upper semi continuity of $\lambda_1$.
\end{proof}

In the case of divisible convex set $\mathcal{C}$, suppose $\Gamma$
divides $\mathcal{C}$, then it suffices to show that
$\Gamma$ is amenable if and only if $\mathcal{C}$ is  amenable.
Hence our results in that case is merely a generalisation of R. Brooks result \cite{brooks} to
this situation.

Hence one gets new examples of convex sets which are not $\delta$-hyperbolic
but have a $\lambda_1>0$: Take the product of a euclidean ball $B_n$ of dimension
$n\geq2$ with a two dimensional triangle or any amenable divisible Hilbert geometry.

Finally let us finish with a question, which is related to the definition
introduced so far

\begin{conj}
  A Hilbert geometry is amenable if and only if it is a $G_n$-polygon.
\end{conj}

The conjecture is trivially true in dimension $2$, thanks to \cite{cvc}.

\section{Proof of the main theorem}
\subsection{Discretisations of Hilbert Geometry} \label{discretisations}

\begin{quote}
In this section we make precise some statements related
to discretisations of Hilbert Geometry, noticeably that they also are of bounded geometry.
\end{quote}

\begin{defi}
  A subset $\mathcal{G}$ of a Hilbert geometry $\mathcal{C}$ is said to be $\varepsilon$-separated,
$\varepsilon>0$, if the distance between any two distinct points of $\mathcal{G}$ is greater than or equal to $\varepsilon$.
\end{defi}

If $\mathcal{G}$ is an $\varepsilon$-separated net, then one always has only a finite numbers of elements
of $\mathcal{G}$ in the ball of radius $r$ centered at a point $x\in \mathcal{C}$, $B(x,r)$.
This is due to the compactness of the balls, which can be covered by a finite number of balls
of radius $\varepsilon/2$. The real difficulty usually lies in obtaining a uniform upper bound
for card $\mathcal{G}\cap B(x,r)$. This is possible thanks to the results in \cite{cvc} and \cite{cvc2}.

\begin{lema}\label{bdgeometry1}
  Let $(\mathcal{C},d_\mathcal{C})$ be a Hilbert geometry and $\mathcal{G}$ an $\varepsilon$-separated subset 
of $\mathcal{C}$, then for all $x \in \mathcal{C}$ and $r>0$
$$
\text{card} \bigl\{\mathcal{G}\cap B(x,r) \bigr\} \leq e^{n\varepsilon}2^n\left(\dfrac{(e^{8r+2\varepsilon}-1)\cdot(e^{\varepsilon+2}-1)}{e^{\varepsilon}-1}\right)^n
$$
\end{lema}

\begin{proof}
  From theorem 9 in \cite{cvc} we have for any hilbert geometry, denoting by $\mu_\mathcal{C}$ its hilbert
measure
$$
\frac{\omega_n}{4^ne^{2nr}}\bigl(\dfrac{e^{2r}-1}{e^{2(r+1)}-1} \bigr)^n \leq \mu_\mathcal{C}\bigl(B(x,r)\bigr)\leq \left(\dfrac{e^{4r}-1}{2}\right)^n \omega_n\text{.}
$$
hence the lemma.
\end{proof}

\begin{defi}
  Let $(\mathcal{C},d_\mathcal{C})$ be a Hilbert geometry. A discretisation of $M$ is a
graph \textbf{G}, determined by an $\varepsilon$-separated subset $\mathcal{G}$ of $\mathcal{C}$ for which
there exists $\rho>0$ such that 
$$
\mathcal{C}=\bigcup_{\xi\in\mathcal{G}}B(\xi,\rho)\text{.}
$$ 
Then $\varepsilon$ is called the separation, and $\rho$ the covering, radius of the discretisation.
The Graph structure \textbf{G} is determined by the collection of $\xi$,
$$
N(\xi):=\{\mathcal{G}\cap B(\xi,3\rho)\}\setminus{\xi} 
$$
for each $\xi\in \mathcal{G}$.
\end{defi}

\begin{rmq}
The choice of the graph structure is such that the graph is always connected.
  Lemma \ref{bdgeometry1} implies that the graph $\mathcal{G}$ is of bounded geometry, i.e., the number
of edges at each vertices is uniformly bounded.
\end{rmq}

\begin{prop}
  Let $(\mathcal{C},d_\mathcal{C})$ be a Hilbert Geometry and \textbf{G} a discretisation of $\mathcal{C}$.
Then there exist $a>1$ and $b>0$ for which
$$
a^{-1}d_\mathcal{C}(\xi_1,\xi_2)\leq d_{\textbf{G}}(\xi_1,\xi_2)\leq a d_\mathcal{C}(\xi_1,\xi_2)+b
$$
for all $\xi_1$,$\xi_2$ in $\mathcal{G}$. Thus $\mathcal{C}$ is quasi-isometric to any
of its discretisations, and any two of its discretisations are quasi-isometric.
\end{prop}
\begin{proof}
  Let $\rho$ be the covering radius of the discretization, and consider a path from
$\xi_1$ to $\xi_2$, then it is clear that 
$$
d_\mathcal{C}(\xi_1,\xi_2)\leq\rho d_{\textbf{G}}(\xi_1,\xi_2)\text{.}
$$
Now consider two points $\xi_1$ and $\xi_2$ in $\mathcal{G}$ and a minimising path in $\mathcal{C}$ from
$\xi_1$ to $\xi_2$. Cut this path into pieces of at most $\varepsilon$ length. This gives less than $d_\mathcal{C}(\xi_1,\xi_2)/\varepsilon +1$
points on the path. Now for each points (excepted the extremities) take the point of $\mathcal{G}$
the closest to it. Thanks to the triangle inequality on can see that we built a path in the graph from $\xi_1$ to
$\xi_2$ with a length less than $d_\mathcal{C}(\xi_1,\xi_2)/\varepsilon +1$.
\end{proof}

\begin{prop}
  Let $(\mathcal{C},d_\mathcal{C})$ be a Hilbert Geometry. Then for any discretisation \textbf{G} of $\mathcal{C}$,
  \begin{enumerate}
      \item \textbf{G} has polynomial
volume growth if and only $\mathcal{C}$ has polynomial volume growth;
  \item \textbf{G} has exponential
volume growth if and only $\mathcal{C}$ has exponential volume growth.
  \end{enumerate}
\end{prop}
\begin{proof}
We do the polynomial growth, the exponential growth goes along the same lines.  
Suppose \textbf{G} has polynomial volume growth, which means that there are constants
$a$ et an integer $d$ such that card $\{\eta\mid d_{\textbf{G}}(\xi,\eta)\leq R  \} \leq aR^d$. Now let us consider
a ball $B(\xi,R)$ in $\mathcal{C}$, then it has a volume less than
$$
\text{card} \{\eta\mid d_{\textbf{G}}(\xi,\eta)\leq R  \} \times\left(\dfrac{e^{4\rho}-1}{2}\right)^n \omega_n\leq a'R^d
$$
by theorem 9 in \cite{cvc}. Now suppose that $(\mathcal{C},d_\mathcal{C})$ has polynomial growth,
which means that there is a constant $A$ and an integer $d$ such that
$$
\mu_\mathcal{C}\bigl(B(\xi,R)\bigr) \leq A R^d\text{.}
$$
Then 
\begin{multline}
\text{card} \{\eta\mid d_{\textbf{G}}(\xi,\eta)\leq R  \} \\
\leq \mu_\mathcal{C}\bigl(B(\xi,R)\bigr)\times\frac{4^ne^{2n\varepsilon}}{\omega_n}
\bigl(\dfrac{e^{2(\varepsilon+1)}-1}{e^{2\varepsilon}-1} \bigr)^n \leq A' R^d.
\end{multline}
\end{proof}

\subsection{Local isoperimetric inequality}\label{isolocal}

\begin{quote}
In this section we study the implications of bounded local geometry property
on the volume of balls, spheres and prove a local isoperimetric inequality "\`a la" Buser
in the setting of Hilbert Geometries.
\end{quote}

First let us show that we have a uniform control on the volume of spheres
in the Hilbert geometries.

To do this we use the following lemma whose proof is in \cite{berckver}
\begin{lema}\label{crofton}
  Let $(\mathcal{C},d_\mathcal{C})$ be a Hilbert Geometry in $\R^n$. Consider two convex
sets $A$ and $B$ inside $\mathcal{C}$, such that $A\subset B$. Let us denote by $\nu_{HT}$
the Holmes-Thompson $n-1$ dimensional measure associated to $\mathcal{C}$. Then
$$
\nu_{HT}(\partial A) \leq \nu_{HT}(\partial B)
$$
Furthermore there exists a constant $C(n)$ such that for the Hausdorff measure one has
$$
\nu_\mathcal{C}(\partial A)\leq C(n) \nu_\mathcal{C}(\partial B)\text{.}
$$ 
\end{lema}

\begin{theo}\label{airedesspheres}
  Let $(\mathcal{C},d_\mathcal{C})$ be a Hilbert geometry, then there are two constants
$C_1(r)>0$ and $C_2(r)<\infty$ such that for any point $x$ in $\mathcal{C}$ if $S(x,r)$ denotes
the sphere of radius $r$ centered at $x$, then
$$
C_1(r)\leq \nu_\mathcal{C}\bigl( S(x,r)\bigr)\leq C_2(r).
$$
\end{theo}
\begin{proof}
  Let us suppose that $C_1(r)=0$. This means that for any $\varepsilon$ there is a point $x_\epsilon$
such that $\nu\bigl(S(x_\varepsilon,r)\bigr)\leq \varepsilon/r$, then for any sphere of radius less than
$r$ centered at $x_e$ the same inequality holds, up to a multiplicative constant, 
thanks to lemma \ref{crofton}. 
Now applying the coarea inequality \cite{cvc2} and \cite{berckver}, one
would obtain a ball of measure less than $C'\cdot\varepsilon$. Hence this would contradict theorem 9 in
\cite{cvc}, which states that there is a lower bound on the hilbert measure of balls of
radius $r$.

Let us now suppose that $C_2(r)=\infty$.
This means that for any $M>0$ there is a point
$x_M$, such that $\nu_\mathcal{C}\bigl(S(x_M,r)\bigr)\geq M/r$, then for any sphere of radius bigger
than $r$ centered at $x_M$ the same inequality holds, thanks to lemma \ref{crofton}.
Again by the coarea inequality, the volume of the ball of radius $2r$ centered
at $x_M$ would have a volume bigger that $C''\cdot M$. This again would contradict
the upper bound of theorem 9 in \cite{cvc}.
\end{proof}

One of the key lemmas in \cite{buser} and \cite{kanai} is a local isoperimetric inequality.
We will need such a lemma, so let us state it in our setting

\begin{lema}[local isoperimetric inequality]\label{localii}
  Let $(\mathcal{C},d_\mathcal{C})$ be a Hilbert geometry, $\varepsilon>0$ and $p\in X$. If $H$ is a smooth hypersurface in the
geodesic ball $B_\varepsilon(p)$ dividing it into two non-empty disjoint domains $D_1$ and $D_2$, then the isoperimetric
inequality
$$
\frac{\nu_\mathcal{C}(H)}{\min\bigl\{\mu_\mathcal{C}(D_1), \mu_\mathcal{C}(D_2)\bigr\}}\geq j(\mathcal{C},\varepsilon)
$$
holds, where $j$ is a positive constant. 
\end{lema}
\begin{proof}
  Let us remark that if $\phi$ is a $C$-lipschitz function from a metric space $(X,d_X)$
to a metric space $(Y,d_Y)$, then denoting by $\mu_{t,X}$ and $\mu_{t,Y}$ their respective $t$-haussdorff measures
one has for any subspace $A$ of $X$, that
$$
\mu_{t,Y}\bigl(\phi(A)\bigr)\leq C^t \mu_{t,X}(A)
$$
Now the Hilbert geometries are of local bounded geometry, hence there is a $C$-bilipshitz function $\varphi$
from $B\varepsilon(p)$ to $\R^n$, thus $\varphi\bigl(B_\varepsilon(p)\bigr)$ is inside the ball of radius $C\varepsilon$ centered
at $\varphi(p)$ and contains the ball of radius $\varepsilon/C$ centered at $\varphi(p)$.
Hence it remains to show that the images of $H$, $D_1$ and $D_2$ satisfy a local
isoperimetric inequality in $\R^n$. But this is the content of the local isoperimetric
inquality of lemma 5.1 in P. Buser's paper \cite{buser}. Now using the fact that
lipschitz hypersurface can be approximated by smooth hypersurfaces one deduces
the local isoperimetric inequality in Hilbert Geometry. 
\end{proof}

\subsection{Discretisations and Isoperimetry}\label{isodiscret}

\begin{quote}
In this section, we show that the positivity of the Cheeger constant of a Hilbert Geometry
is the same as the positivity of the Cheeger constants of its dicretisations. The results
follows from the bounded geometry of Hilbert geometries. 
This is quite standard in the setting
of Riemannian geometry and the proof is similar.
\end{quote}

First we must recall what we call Cheeger constant of a graph
\begin{defi}
  The \textit{cheeger constant} of a graph $\textbf{G}$ is
$$
I_\infty(\textbf{G})=\inf\bigl\{  \frac{|\partial F|}{|F|} \big| F \text{is finite and non-empty subset of vertices of} \textbf{G}  \bigr\} 
$$
where $\partial F$ denotes the set of points at a distance less than one from a point of $F$,
and which are not in $F$. As usual we denote by $|F|$ the cardinal of $F$.
\end{defi}

Now let us state the main result of this section
\begin{theo}\label{cheeger}
  Let $(\mathcal{C},d_{C})$ be a Hilbert geometry. Then its Cheeger constant is positive
if and only if the Cheeger constant of any discretisation is positive. 
\end{theo}

This theorem must be linked with
the results of M. Kanai \cite{kanai} related to Riemannian manifolds with Ricci
curvature bounded from below and positive injectivity radius.


\begin{proof}
  Suppose that $I_\infty(\mathcal{C})>0$. We may work with any discretisation. Let us consider
a discretisation \textbf{G} with separation constant $\varepsilon>0$ and covering radius $\rho=R$. To show that
$I_\infty(\textbf{G})>0$ it suffices to prove the existence of positive constants $C_1$ and $C_2$ such that given
any $\mathcal{K}\subset \mathcal{G}$ we may find $\Omega\subset \mathcal{C}$ 
for which
\begin{equation}
  \label{eqcheeger1}
  \nu_\mathcal{C}(\partial\Omega)\leq C_1 \text{card} \partial \mathcal{K}\text{,}
\end{equation}
and
\begin{equation}
  \label{eqcheeger2}
  \mu_\mathcal{C}(\Omega)\geq C_2 \text{card} \mathcal{K}\text{.}
\end{equation}

Given a finite subset $\mathcal{K}$, set
$$
\Omega:=\bigcup_{\xi\in \mathcal{K}} B(\xi,R)\text{.}
$$
Let $M(\varepsilon,R)$ be an upper bound on the ratio of the volume of a disk of  radius $R$
by the volume of a disk of radius $\varepsilon$ (This is also an upper bound of the maximum number
of $\varepsilon$-separated points in a disk of radius $R$).
This bound exists thanks to lemma \ref{bdgeometry1}. Let us also denote by $V_R$ the infimum
of the volume of a ball of radius $R$ in $\mathcal{C}$ (which is note zero thanks to theorem 9 in \cite{cvc}).
Then
\begin{multline}
\sum_{\xi\in \mathcal{K}} \mu_C\bigl(B(\xi,R)\bigr)\leq M(\varepsilon,R)\sum_{\xi\in \mathcal{K}} \mu_C\bigl(B(\xi,\varepsilon)\bigr)= M(\varepsilon,R)\mu_C\Bigl(\bigcup_{\xi\in \mathcal{K}} B(\xi,\varepsilon)\Bigr)\\
\leq M(\varepsilon,R) \mu_C\Bigl(\bigcup_{\xi\in \mathcal{K}} B(\xi,R)\Bigr) = M(\varepsilon,R) \mu_C(\Omega)\text{.}
\end{multline}
thus we obtain
$$
V_R\text{card} \mathcal{K} \leq M(\varepsilon,R) \mu_C(\Omega)
$$
which corresponds to \ref{eqcheeger2}.
For the upper bound on $\nu(\partial\Omega)$, we claim that
$$
\partial\Omega \subset \bigcup_{\xi\in \partial(\mathcal{G}\setminus\mathcal{K})} S(\xi,R)\text{.}
$$
To see this, remark that if $x\in \partial\Omega$, then $d_\mathcal{C}(x,\xi)\geq R$ for all $\xi\in \mathcal{K}$, and there exists
$\xi_0\in \mathcal{K}$ such that $x\in S(\xi_0,R)$ (for if one of these conditions fails, $x$ is either outside
or inside $\Omega$). But by definition, there must exist $\xi'\in \mathcal{G}$ such that $d_\mathcal{C}(x,\xi')< R$,
which implies $\xi'\not\in \mathcal{K}$. However $d_\mathcal{C}(\xi_0,\xi')<2R$, which implies that $\xi_0\in N(\xi')$
in other words $\xi_0\in \partial(\mathcal{C}\setminus\mathcal{K})$.

Therefore using \ref{airedesspheres}, and letting $m$ be the maximum number
of points in the neighbourhood of a point in $\mathcal{G}$, we get
$$
A(\partial\Omega)\leq C_2(R) \text{card} \partial\bigr(\mathcal{G}\setminus\mathcal{K}) \leq m C_2(R) \text{card} \mathcal{K}
$$

Assume now that $I_\infty(\textbf{G})>0$. Suppose we are given $\Omega$, with compact closure and
$C^\infty$ boundary in $\mathcal{C}$. Set
\begin{eqnarray*}
  \mathcal{K}_0&:=& \Bigl\{\xi\in \mathcal{G}\mid  \mu_\mathcal{C}\bigl(\Omega\cap B(\xi,\rho)\bigr)> \mu_\mathcal{C}\bigl(B(\xi,\rho) \bigr) \Bigr\} \\
  \mathcal{K}_1&:=&\Bigl\{\xi\in \mathcal{G}\mid  \mu_\mathcal{C}\bigl(\Omega\cap B(\xi,\rho)\bigr)\leq \mu_\mathcal{C}\bigl(B(\xi,\rho) \bigr) \Bigr\} 
\end{eqnarray*}
Both $\mathcal{K}_0$ and $\mathcal{K}_1$ are contained in $\Omega_\rho$, the set
of points at distance less or equal to $\rho$ from $\Omega$. Furthermore for at least
 one of $j=0,1$ we have
\begin{equation}\label{eqcinqtroissix}
\dfrac{\mu_C(\Omega)}{2}\leq \mu_C\biggl(\Omega\cap \bigcup_{\xi\in \mathcal{K}_j} B(\xi,\rho) \biggr) \text{.}
\end{equation}

Assume equation (\ref{eqcinqtroissix}) is valid for $j=0$.
Denote by $V_\mathcal{C}(\rho)$ the upper bound on the volume of balls of
radius $\rho$ in $\mathcal{C}$. First notice that
$$
\dfrac{\mu_C(\Omega)}{2} \leq \sum_{\eta\in \mathcal{K}_0} \mu_C\bigl((\Omega\cap B(\eta,\rho)\bigr) \leq V_\mathcal{C}(\rho)\text{card }\mathcal{K}_0 
$$
thus it suffices to give a lower bound of $\nu(\partial\Omega)$ by a multiple of card $\mathcal{K}_0$,
the multiple being, of course independent of $\mathcal{K}_0$. To do this
define $H\subset\Omega_\rho$ by
$$
H:=\Bigl\{x\in M\mid \mu_\mathcal{C}\bigl(B(x,\rho)\bigr)/2=\mu_\mathcal{C}\bigl(\Omega\cap B(x,\rho)\bigr) \Bigr\}\text{.}
$$

For each $\xi\in \partial\mathcal{K}_0$ there exists $\eta\in N(\xi)$, $\eta\in \mathcal{K}_0$. From the defintion of $N(\xi)$
it follows that
$$
d_\mathcal{C}(\xi,\eta)< 3\rho\text{.}
$$

The definitions of $\mathcal{K}_0$ and $\partial\mathcal{K}_0$ imply
$$
\mu_\mathcal{C}\bigl(\Omega\cap B(\eta,\rho)\bigr)> \mu_\mathcal{C}\bigl(B(\eta,\rho) \bigr), \quad %
\mu_\mathcal{C}\bigl(\Omega\cap B(\xi,\rho)\bigr)\leq \mu_\mathcal{C}\bigl(B(\xi,\rho) \bigr)
$$

Thus by continuity of the volume, the line between $\xi$ and $\eta$ contains an element
$\zeta\in H$, which implies $\partial\mathcal{K}_0\subset H_{3\rho}$, and
$$
\bigcup_{\xi\in \partial\mathcal{K}_0} B(\xi,\rho) \subset H_{4\rho}
$$

Now let $Q$ be a maximal $2\rho$-separated subset of $H$ thus
$$
\bigcup_{\xi\in \partial\mathcal{K}_0} B(\xi,\rho) \subset Q_{6\rho}
$$
which implies
\begin{eqnarray*}
  V_\rho\text{card }\partial\mathcal{K}_0 &\leq& \sum_{\xi\in \partial\mathcal{K}_0} \mu_\mathcal{C}\bigl(B(\xi,\rho)\bigr)\\
&\leq & M_{\epsilon,\rho} \sum_{\zeta\in Q}  \mu_\mathcal{C}\bigl(B(\zeta,6\rho)\bigr)\\
\text{by theorem 9 in \cite{cvc}}&\leq & M_{\epsilon,\rho} const. \sum_{\zeta\in Q}  \mu_\mathcal{C}\bigl(B(\zeta,\rho)\bigr)\\
&=& 2M_{\epsilon,\rho} const. \sum_{\zeta\in Q} \mu_\mathcal{C}\bigl(\Omega\cap B(\zeta,\rho)\bigr)\\
\text{by lemma \ref{localii}} &\leq& 2M_{\epsilon,\rho} const.' \sum_{\zeta\in Q} \nu_\mathcal{C}\bigl(\partial\Omega\cap B(\zeta,\rho)\bigr)\\
&\leq& 2M_{\epsilon,\rho}^2 const.' \nu_\mathcal{C}(\partial\Omega)\text{.}
\end{eqnarray*}

Now assume equation (\ref{eqcinqtroissix}) is valid for $j=1$. Then we have from
lemma \ref{localii}
\begin{multline}
  \label{cas1}
  \dfrac{\mu_C(\Omega)}{2} \leq \sum_{\xi\in \mathcal{K}_1} \mu_\mathcal{C}\bigl( \Omega\cap B(\xi,\rho)\bigr)\\
   \leq \text{const.} \sum_{\xi\in \mathcal{K}_1} \nu_\mathcal{C}\bigl( \partial\Omega\cap B(\xi,\rho)\bigr) \leq \text{const.} M_{\varepsilon,\rho} \sum_{\xi\in \mathcal{K}_1} \nu_\mathcal{C}\bigl( \Omega\cap B(\xi,\varepsilon)\bigr)\\
   = \text{const.} M_{\varepsilon,\rho} \nu_\mathcal{C}\bigl(\Omega\cap \bigcup_{\xi\in \mathcal{K}_1} B(\xi,\varepsilon) \bigr)\leq \text{const.} M_{\varepsilon,\rho} \nu_\mathcal{C}(\partial\Omega)
\end{multline}
which finishes the proof.
\end{proof}

\subsection{Isoperimetry and bottom of spectrum}\label{isobot}

\begin{quote}
In this section we  recall how the Cheeger constant of a discretisation is related to its spectral radius and to amenability.
\end{quote}

To go further into the subject one should consult \cite{cgh}, we will extract from this paper
the notion needed here.

Let us first start by recalling that on a locally finite graph  \textbf{G}, whose set of  vertices is $\mathcal{G}$, there is a natural
\textsl{simple random walk} with corresponding \textsl{Markov operator T}. We will also suppose that  \textbf{G} is connected and of bounded degree, which is the case for our discretisations as we saw in the previous
sections.
Then one can consider the Hilbert space $l^2(\mathcal{G},deg)$ of functions $h$ from the vertices
$\mathcal{G}$ to $\C$ such that $\sum_{x\in \mathcal{G}} deg(x)|h(x)|^2<\infty$, and the bounded self-adjoint operator
$T$ defined on this Hilbert space by
$$
(Th)(x)=\frac{1}{deg(x)}\sum_{y\sim x}h(y)
$$
where $y\sim x$ indicates a summation over the neighbours $y\in N(x)$ of the vertex $x$.
The \textsl{spectral radius} of \textbf{G} is
\begin{eqnarray*}
  \rho(\textbf{G})&=&\sup\bigl\{ \langle h,Th\rangle\big| h\in l^2(X), ||h||_2\leq 1 \bigr\}\\ 
&=&\sup\bigl\{ |\lambda| \mid \lambda \text{ is in the spectrum of } T\bigr\}\text{.} 
\end{eqnarray*}

With this notions in mind one must also notice that $1-T$ is a natural analogue on \textbf{G} of
a Laplacian, so that $1-\rho(\textbf{G})$ is usually referred to as the bottom of its spectrum.

\begin{rmq}\label{fpropre}
It is also known that, for a any real number $\lambda\geq\rho(\textbf{G})$ there exists $F\colon \mathcal{G}\to\mathopen]0,\infty\mathclose[$ such that
\begin{equation}
\frac{1}{deg(x)}\sum_{y\sim x}F(y)=\lambda F(x)\text{.}
\end{equation}
Actually this is an equivalence.
\end{rmq}

Another equivalent definition of $\rho(\textbf{G})$ is the following. For $x,y\in \mathcal{G}$
and for any integer $n\leq0$, denote by $p^{(n)}(x,y)$ the probability that a simple random walk starting
at $x$ is at $y$ after $n$ steps. Then one has also
$$
\rho(\textbf{G})= \limsup_{n\to \infty }\sqrt[n]{p^{(n)}(x,y)}
$$

To conclude our paper it remains to finish the exploration of the 
link between the bottom of the spectrum and the cheeger constant.
This is the content of the following two results, stated without proof 
(see \cite{cgh} and references werein).

This first lemma is a kind of inverse Cheeger inequality.
\begin{lema}
  For a graph \textbf{G} which is regular of degree $d\geq2$, one has
$$
I_\infty(\textbf{G})\geq 4\frac{1-\rho(x)}{\rho(x)}
$$
\end{lema}

Finally the missing piece of our puzzle is the following one
\begin{theo}\label{link}
  Let \textbf{G} be a connected graph of bounded degree. The following are equivalent
  \begin{enumerate}
      \item \textbf{G} is not amenable;
      \item $I_\infty(\textbf{G})>0$;
      \item $\rho(\textbf{G})<1$;
      \item  $p^{(n)}(x,y) = o(s^n)$ for some $s\in \mathopen]0,1 \mathclose[$ and for all $x,y \in \mathcal{G}$.
  \end{enumerate}
Should one of this be true, then the simple random walk on \textbf{G} is transient.
\end{theo}

\subsection{Bottom of the spectrum and discretisations}\label{lambdadiscret}

\begin{quote}
In this section, we show that the positivity of the bottom of spectrum of a Hilbert Geometry
is the same as the positivity of the bottom of the spectrum of its dicretisations. Once
again the path is standard in Riemannian geometry, and follows by the bounded geometry property.
\end{quote}

Now let us state the main result of this section
\begin{theo}\label{cheeglap}
  Let $(\mathcal{C},d_{C})$ be a Hilbert geometry. The bottom of the spectrum
of $\mathcal{C}$ is positif if and only if the spectral radius of any discratisation
is stritly smaller than $1$.
\end{theo}

To prove this theorem we will need to raise functions on the discretisations
to functions on the convex. We do as follows.

Consider $(\phi_\xi)_{\xi\in \mathcal{G}}$ a partition of unity on $\mathcal{C}$ subordinate to the locally finite
cover $\{B(\xi,2\rho)\}_{\xi\in \mathcal{G}}$, and such that $\phi_\xi=1$ on $B(\xi,\rho)$.
Then for each $f\colon \mathcal{G}\to \R$ we define its smoothing
$F=\mathcal{S}f\colon \mathcal{C}\to \R$ by
$$
(Sf)(x)=\sum_{\xi\in \mathcal{G}} \phi_\xi(x)f(\xi)\text{.}
$$ 

our main claim, whose proof we postpone, is the following
\begin{lema}[smoothing lemma]\label{smoothing}  
Let $(\mathcal{C},d_{C})$ be a Hilbert geometry and $\mathcal{G}$ one discretisation
of $\mathcal{C}$. Let $\mathcal{S}$ be the smoothing operator defined as above
and $T$ the markov operator associated to the simple random walk on $\mathcal{G}$.
There exists two constants $C_1$ and $C_2$ such that
\begin{eqnarray}
  ||f||_2^2&\leq& C_1 ||\mathcal{S}f||_2^2 \\
  ||d\mathcal{S}f||_2^2 & \leq & C_2 \langle(1-T)f,f\rangle
\end{eqnarray}
\end{lema}

\begin{proof}[Proof of theorem \ref{cheeglap}]
  Suppose $\rho(\textbf{G})<1$, then by theorem \ref{link}  the Cheeger constant
of the graph is positive, 
$I_\infty(\textbf{G})>0$. Now theorem \ref{cheeger} implies that it is also the
case for the cheeger constant of our hilbert geometry $I_\infty(\mathcal{C})$. 
Finally using the Cheeger inequality
proved in \cite{cvc} we obtain that $\lambda_1(\mathcal{C})>0$.

Assume now that $\rho(\textbf{G})=1$. Hence for any $\lambda\geq1$, by the remark
\ref{fpropre} there exists a function $F\colon\mathcal{G}\to\R^+_*$ such that
$$
\frac{1}{deg(x)}\sum_{y\sim x}F(y)=\lambda F(x)
$$ 
(As our discretisations are of bounded degree, without loss we can consider
that $deg(x)$ is a constant, and take this constant equal to $1$.)

We can rewrite this last equality under the following form
$$
\langle(1-T)F,F\rangle=(1-\lambda)||F||^2
$$
Hence by taking cut off functions and $\lambda=1$ we deduce the existence of a familly of functions $f_n$
with compact support on $\mathcal{G}$, such that
$$
\frac{\langle(1-T)f_n,f_n\rangle}{||f_n||^2}\leq \frac{1}{n}\text{.}
$$
now we can easily conclude thanks to the smoothing lemma \ref{smoothing} that $\lambda_1(\mathcal{C})=0$.
\end{proof}

\begin{proof}[Smoothing lemma's proof]
  Recall that there is a constant $V_\rho$ wich is a lower bound on the volume of balls of radius
$\rho$ in $\mathcal{C}$, thanks to theorem 9 in \cite{cvc}. Hence we have
\begin{eqnarray*}
  \int_\mathcal{C} (\mathcal{S}f)^2 d\mu_\mathcal{C}(x) &\geq& \sum_\xi \int_{B(\xi,\epsilon/2)} (\mathcal{S}f)^2 d\mu_\mathcal{C}(x)\\
&\geq&\sum_\xi \int_{B(\xi,\epsilon/2)}\phi_\xi^2(x)f^2(\xi)d\mu_\mathcal{C}(x)\\
&=& \sum_\xi \int_{B(\xi,\epsilon/2)}f^2(\xi)d\mu_\mathcal{C}(x)\geq V_\epsilon  \sum_\xi f^2(\xi)\text{.}
\end{eqnarray*}

Now let us consider the differentials and $V\in \R^n$
$$
d (\mathcal{S}f)_x\cdot V = \sum_\xi f(\xi) d(\phi_\xi)_x\cdot V =  \sum_{\xi\in B(x,2\rho)}f(\xi) d(\phi_\xi)_x\cdot V
$$
Given $x$ there exists $\eta_x \in \mathcal{G}\cap B(x,\rho)$ hence
\begin{equation}\label{unvoisinage}
  d (\mathcal{S}f)_x\cdot V =\sum_{\xi\in B(\eta_x,3\rho)}f(\xi) d(\phi_\xi)_x\cdot V
\end{equation}
and since  the $\phi_\xi$ are a partition of unity we have $\sum_{\xi\in B(\eta_x,3\rho)} \phi_\xi(x)=1$,
and differentiating one gets $\sum_{\xi\in B(\eta_x,3\rho)} d(\phi_\xi)_x=0$.  Using this in (\ref{unvoisinage}) we
finally obtain
\begin{equation*}
 d (\mathcal{S}f)_x\cdot V= \sum_{\xi\in B(\eta_x,3\rho)}\bigl(f(\xi)-f(\eta_x)\bigr) d(\phi_\xi)_x\cdot V\text{.}
\end{equation*}
Therefore
$$
F_\mathcal{C}^*(x,d(\mathcal{S}f)_x)\leq C \sum_{\xi\in B(\eta_x,3\rho)}\bigl|f(\xi)-f(\eta_x)\bigr|
$$
which implies for any $x\in B(\eta,\rho)$, $\eta\in \mathcal{G}$
that
$$
(F_\mathcal{C}^*)^2(x,d(\mathcal{S}f)_x)\leq C' \sum_{\xi\in B(\eta,3\rho)}\bigl|f(\xi)-f(\eta)\bigr|^2= C''|df|^2(\eta)
$$
and now using the fact that the hilbert geometry is quasi-isometric to its
discretisation we deduce the inequality which follows
$$
\int_{B(0,R)} (F_\mathcal{C}^*)^2(x,d(\mathcal{S}f)_x) d\mu_\mathcal{C}(x) \leq C_2 \int_{\beta(\eta_0,R+1)} |df|^2dV
$$
and taking $R\to \infty$ we finally obtain
$$
||d\mathcal{S}f||_2^2  \leq  C_2 \int_{\beta(\eta_0,R+1)} |df|^2dV= C_2 \langle(1-T)f,f\rangle\text{.}
$$
\end{proof}

\nocite{*}
\bibliographystyle{amsalpha}
\bibliography{enthilb}

\end{document}